\documentclass[a4paper,twoside]{article}

\newcommand{\heute}{March 22, 2001}
\usepackage{amsmath,amsthm,amsfonts,latexsym,amscd}
\usepackage{amsthm,amssymb}

\newcommand{\href}[2]{#1}

\swapnumbers
\theoremstyle{plain}
\newtheorem{theorem}{Theorem}[section]
\newtheorem{lemma}[theorem]{Lemma}
\newtheorem{corollary}[theorem]{Corollary}
\newtheorem{proposition}[theorem]{Proposition}

\theoremstyle{definition}

\newtheorem{example}[theorem]{Example}

\newtheorem{remarks}[theorem]{Remarks}

\newcommand{\NeumannN}{\mathcal{N}}

\newcommand{\reals}{\mathbb{R}}
\newcommand{\complexs}{\mathbb{C}}
\newcommand{\naturals}{\mathbb{N}}
\newcommand{\integers}{\mathbb{Z}}
\newcommand{\rationals}{\mathbb{Q}}

\newcommand{\abs}[1]{\left\lvert#1\right\rvert} %absolute value

\newcommand{\semiProd}{\rtimes}

\DeclareMathOperator{\de}{d}

\DeclareMathOperator{\tr}{tr}
\DeclareMathOperator{\pr}{pr}
\DeclareMathOperator{\avg}{avg}  %average
\newcommand{\forget}[1]{}

\newcommand{\innerprod}[1]{\langle #1 \rangle}

{\catcode`@=11\global\let\c@equation=\c@theorem}

\allowdisplaybreaks[2]

\begin{document}

\pagestyle{myheadings}
\markboth{Spectral measure and wreath products}{Warren Dicks and Thomas Schick}

\title{The spectral measure of certain elements of the complex group
  ring of a wreath
  product}
\author{
Thomas Schick\thanks{
e-mail: thomas.schick@math.uni-muenster.de
\protect\\
www:~http://wwwmath.uni-muenster.de/u/lueck/
\protect\\
Fax: ++49 -251/83 38370}
}
\author{Warren Dicks\thanks{
e-mail: dicks@mat.uab.es
\protect\\
www:~http://mat.uab.es/dicks/
\protect\\
Fax: ++34 -93/581 2790\protect\\
Research funded by the DGI (Spain) through Grant BFM2000-0354.
}\\
Departament de Matem\`atiques\\
Universitat Aut\`onoma de Barcelona\\
08193 Bellaterra (Barcelona)\\
Spain
\and
Thomas Schick\thanks{
e-mail: thomas.schick@math.uni-muenster.de
\protect\\
www:~http://wwwmath.uni-muenster.de/u/lueck/
\protect\\
Fax: ++49 -251/83 38370}\\
%Penn State University and Universit\"at M\"unster
%FB Mathematik --- Universit\"at M\"unster\\
%Einsteinstr.~62 --- 48159 M\"unster\\
%Germany
FB Mathematik\\
Universit\"at M\"unster\\
Einsteinstr.~62\\
48159 M\"unster\\
Germany
}

\date{Last edited: \heute.  Last compiled: \today.}

\maketitle

\begin{abstract}
  We use elementary methods to compute the $L^2$-dimension of the
  eigen\-spaces of the Markov operator on the lamplighter group and
  of generalizations of this operator on other groups. In particular,
  we give a transparent explanation of the spectral measure of the
  Markov operator on the  lamplighter group found by Grigorchuk-Zuk
  \cite{Grigorchuk-Zuk(2000)}. The latter result was used by
  Grigorchuk-Linnell-Schick-Zuk
  \cite{Grigorchuk-Linnell-Schick-Zuk(2000)} to produce a
  counterexample to a strong version of the Atiyah conjecture about
  the range of $L^2$-Betti numbers.

  We use our results to construct manifolds with certain $L^2$-Betti
  numbers (given as convergent infinite sums of rational numbers) which
  are not obviously rational, but we have been unable to determine whether
  any of them are irrational.
\end{abstract}

\section{Notation and statement of main result}
\label{sec:notat-stat-main}

In this section we introduce notation that will be fixed throughout and will be
used in the statement of the main result.

Let $U$ denote a discrete group with torsion.

Let $e$ be a nontrivial projection (so $e=e^*=e^2$, $e \ne 0,1$) in
$\complexs[U]$.
For example, $U$ could be finite and nontrivial,  and $e$ could be the
`average' of
the elements of $U$,
\begin{equation*}
\avg(U):= \frac{1}{\abs{U}}\sum_{u\in U} u.
\end{equation*}
This will be the example we shall make the most use of.

Let $W = W(U,e)$ denote the inverse of the coefficient of 1 in the
expression of $e$
as a $\complexs$-linear combination of elements of $U$.  By results of
Kaplansky and
Zaleskii, $W$ is a rational number greater than 1.   For example, if
$U$ is finite and nontrivial, and $e=\avg(U)$,  then $W=\abs{U}$.

For integers $m$, $n$, with $1 \le m \le n-1$, let
$\lambda_{m,n}:= 2 \cos(\frac m n \pi)$.

For any integer $n \ge 2$, let
$M_n:=\{\lambda_{m,n}\mid 1 \le m \le n-1, m\text{ coprime to } n\}$.

We write
\begin{equation*}
U \wr \integers := (\oplus_{i\in\integers} U)\semiProd C_\infty,
\end{equation*}
where $C_\infty$ denotes an infinite cyclic group with generator $t = t_U$
which acts
on $\oplus_{i\in\integers} U$ by the shift,
i.e.~$t^{-1}((g_n)_{n\in\integers})t=(g_{n-1})_{n\in\integers}$. For each
$u \in U$,
let $a_u$ denote $(\dots,1,u,1,\dots)\in
\oplus_{i\in\integers}U$ where $u$ occurs with index $0$.  Throughout, we
identify
$u$ with $a_u$.  Thus $U$ is a subgroup of $U \wr \integers$. Notice that  $U\wr
\integers$ is generated by $t$ and $U$.

Set
\begin{equation*}
  T = T(U, e) := \left(et + t^{-1}e\right) \in \complexs [U \wr \integers].
\end{equation*}
If $U$ is finite and nontrivial, and $e=\avg(U)$,  then $T$ is two times
the Markov operator of $U\wr \integers$ with respect to the
symmetric set of generators $\{ut,(ut)^{-1}\mid u\in U\}$.

Let $\NeumannN(U\wr \integers)$ denote the (von Neumann) algebra of bounded
linear operators on the Hilbert space  $l^2(U\wr \integers)$  which commute with
right multiplication by elements of $U\wr \integers$. We identify each element
$x$ of $\complexs[U\wr \integers]$ with an element of $l^2(U\wr \integers)$
in the
natural way, and also with the element of $\NeumannN(U\wr \integers)$
given by left multiplication by $x$.  Thus $\complexs[U\wr \integers]$ is
viewed as a
subset of $l^2(U\wr \integers)$ and as a subalgebra of $\NeumannN(U\wr
\integers)$.
For $a \in \NeumannN(U\wr
\integers)$ the (regularized) trace of $a$ is defined as
\begin{equation*}
  \tr_{U\wr \integers}(a):= \innerprod{a(1),1}_{l^2(U\wr \integers)}.
\end{equation*}

  Similar notation applies for any group.

  Note that, if $a\in \NeumannN(U\wr \integers)$ leaves invariant
  $l^2(G)$ for a subgroup $G$, then we can consider $a$ to be an element of
  $\NeumannN(G)$, and here $\tr_G(a)$ and $\tr_{U\wr \integers}(a)$ coincide.

  Note also that, if $a$ lies in $\complexs[U\wr \integers]$, then
  $\tr_{U\wr \integers}(a)$ is the coefficient of 1 in the expression of $a$
  as a $\complexs$-linear combination of elements of $U\wr \integers$.

The element (left multiplication by) $T$ of $\NeumannN(U\wr\integers)$ is
self-adjoint.  For each $\mu \in \reals$, let
$\pr_\mu \colon l^2(U\wr \integers) \to l^2(U\wr \integers)$ denote the
orthogonal
projection onto $\ker (T -\mu)$, so $\pr_\mu \in \NeumannN(U\wr
\integers)$.  The
number
\[
\dim_{U\wr \integers}\ker( T-\mu) := \langle \pr_\mu(1), 1 \rangle _{l^2(U\wr
\integers)}=
\tr_{U\wr \integers}(\pr_{\mu})
\]
is called the $L^2$-multiplicity of $\mu$ as an eigenvalue of $T$.

Our main result is the following.
\begin{theorem}\label{bigmain}
  With all the above notation, for any $\mu \in \reals$,
\begin{equation*}
    \dim_{U \wr \integers} \ker(T-\mu) =\begin{cases}
      \frac{(W-1)^2}{W^{n}-1} & \text{if }n \ge 2 \text{ and }\mu\in M_n,\\
  0 &\text{if }\mu\notin \bigcup_{n\ge 2} M_n.
\end{cases}
\end{equation*}
Moreover, $l^2(U\wr \integers)$ is the Hilbert sum of the eigenspaces of
$T$, i.e.~the
spectral measure of $T$ off its eigenspaces is zero.
\end{theorem}

In \cite[Corollary 3]{Grigorchuk-Zuk(2000)}, Grigorchuk-Zuk proved the case
of this
result in which  $U$ is (cyclic) of order two and $e=\avg(U)$, so $W = 2$.
This was
used in \cite{Grigorchuk-Linnell-Schick-Zuk(2000)} to give a counterexample to a
strong version of the Atiyah conjecture about the range of $L^2$-Betti
numbers.
The argument in \cite{Grigorchuk-Zuk(2000)} is based on automata and actions on
binary trees, while our proof is based on calculating traces of projections
in the
group ring $\complexs[U\wr \integers]$.

\section{Preliminary matrix calculations}

In this section, we introduce more notation which will be used throughout,
and verify
some identities which will be used in the proof.

For positive integers $i$, $j$, let
\begin{equation*}
\alpha_{i,j} := \delta_{\abs {i-j},1} =
\begin{cases}
  1 & \text{if } i -j = \pm 1,\\
  0 & \text{otherwise.}
\end{cases}
\end{equation*}
For each integer $n \ge 2$, let $A_n$ denote the $n-1 \times n-1$ matrix
\begin{equation*}
A_n = (\alpha_{i,j})_{1 \le i,j\le n-1} =
 \begin{pmatrix}
      0 & 1 & 0 &\hdotsfor{4}\\
      1 & 0 & 1 & 0 & \hdotsfor{3}\\
      0 & 1 & 0 & 1 & 0 & \hdotsfor{2}\\
      \hdotsfor[2]{7}\\
      \hdotsfor{3} & 0 & 1 & 0 & 1\\
      \hdotsfor{4} & 0 & 1 & 0
    \end{pmatrix}.
\end{equation*}

Recall that $\lambda_{m,n}$ denotes $2 \cos(\frac{m}{n}\pi)$.

\begin{lemma}
    For each $n \ge 2$,  the family of eigenvalues of $A_n$, with
multiplicities,
    is $\{\lambda_{m,n} \mid 1 \le m \le n-1\}$.
  \end{lemma}
  \begin{proof}
    For a complex number $\mu$ different from $0, 1, -1$, one checks
    immediately
    by induction on $n$, and
    determinant expansion of the first row,  that
    \begin{equation*}
      \det(A_n+(\mu+\mu^{-1})I_{n-1}) = \frac{\mu^n - \mu^{-n}}{\mu -
        \mu^{-1}}.
    \end{equation*}
    Now, for $1 \le m \le n-1$, taking  $\mu = -e^{\frac{m}{2n} 2\pi i}$
    shows that  $\lambda_{m,n}$ is an eigenvalue of $A_n$.  Since we have $n-1$
    distinct eigenvalues for $A_n$, they all have multiplicity one.
\end{proof}

For $n \ge 2$,  $A_n$ is a real symmetric matrix, so there exists a real
orthogonal matrix $B_n = (\beta_{i,j}^{(n)})_{1 \le i,j \le n-1}$ such that
$B_nA_nB_n^{*}$ is a diagonal matrix $D_n$; here the diagonal entries are
$\lambda_{m,n}$, $1 \le m \le n-1$, and we may assume the entries occur in
this order,
so $D_n = (\delta_{i,j}\lambda_{j,n})_{1 \le i,j \le n-1}$.
Since $B_nB_n^* = I_{n-1}$ and $B_nA_n = D_nB_n$ we have the identities
\begin{equation}\label{eq:orthogonal}
\sum_{j = 1}^{n-1} \beta_{i,j}^{(n)} \beta_{k,j}^{(n)} =
\delta_{i,k}, \quad 1 \le i,k \le n-1,
\end{equation}
\begin{equation}\label{eq:diagonal}
 \sum_{j = 1}^{n-1} \beta_{i,j}^{(n)}\alpha_{j,k}
= \lambda_{i,n} \beta_{i,k}^{(n)}, \quad 1 \le i,k\le n-1.
\end{equation}

\section{Proof of the main result}

We shall frequently use the following, which is well known and easy to prove.

\begin{lemma}\label{lem:tensor_product_trace}
  Let $G$ and $H$ be discrete groups, and let  $p\in\NeumannN(G)$ and
  $q\in\NeumannN(H)$. Embed $G$ and $H$ in
  the canonical way into $G\times H$, so $p$ and $q$ become
  elements of $
  \NeumannN(G\times H)$. Then
  \begin{equation*}
    \tr_{G\times H}(pq)=\tr_G(p)\cdot \tr_H(q).\qed
  \end{equation*}
\end{lemma}

We need even more notation.

  For each $i \in \integers$, we define, in $\complexs[U\wr \integers]$,
$e_i:= t^{-i}et^i$ and $f_i:=1-e_i$.

It is easy to see that all the $e_i$, $f_j$ are projections which commute
with each
other;  moreover,
\begin{equation}\label{eq:first_traces}
\tr_{U \wr \integers} (e_i) = \tr_{U \wr \integers} (e) = \frac 1 W \text{
and  }
\tr_{U \wr \integers} (f_i)= 1 - \frac 1 W.
\end{equation}

For $n \ge 2$, let $q_n := f_1 e_2 e_3 \cdots e_{n-2} e_{n-1} f_{n}.$  It
is clear
that $q_n$ is a projection.  Moreover, the  factors lie in
$\complexs[t^{-i}Ut^i]$,
$1\le i \le n$, so, by Lemma~\ref{lem:tensor_product_trace},
\begin{equation*}
\tr_{U \wr \integers}(q_n) =
\tr_{U \wr \integers}(f_{1}) \tr_{U \wr \integers}(e_{2}) \cdots
\tr_{U \wr \integers} (e_{n-1})\tr_{U \wr \integers} (f_{n}).
\end{equation*}
By (\ref{eq:first_traces}),
\begin{equation}\label{eq:trace}
\tr_{U \wr \integers}(q_n) =  (1 - \frac 1 W)^2(\frac{1}{W})^{n-2} =
\frac{(W-1)^2}{W^n}.
\end{equation}

\begin{lemma}\label{lem:orthogonality}
If $1 \le m < n$ and $1 \le m' < n'$ then $$q_{n'}t^{-m'}t^mq_n =
\delta_{n,n'}\delta_{m,m'}q_n.$$
\end{lemma}

\begin{proof} Note that $t^m q_n t^{-m} = f_{1-m} e_{2-m} \cdots e_{n-m-1}
f_{n-m}$,
and this is a projection. Thus
\begin{equation*}
(t^m q_n t^{-m} \mid n \ge 2, 1 \le m < n)
= (f_{-i}e_{-i+1}\cdots e_{j-1} f_{j} \mid -i \le 0, 1 \le j).
\end{equation*}
This is a family of pairwise orthogonal projections, since, if $-i,-i' \le
0, 1 \le
j,j'$,  then either $(i,j) = (i',j')$,  or the product of $f_{-i}e_{-i+1}\cdots
e_{j-1} f_{j}$ and $f_{-i'}e_{-i'+1}\cdots e_{j'-1} f_{j'}$ is zero since
it contains
a factor $e_\alpha f_\alpha = 0$ for at least one $\alpha \in \{-i,-i', j,j'\}$.
Since $t$ is invertible, the result follows.
\end{proof}

Notice that, for $1\le m<n$,
\begin{equation*}
\begin{split}
T(t^m q_n) &=  ett^m q_n  + t^{-1}et^m q_n \\
&= t^{m+1}e_{m+1}q_n + t^{m-1}e_{m}q_n \\
&=  t^{m+1}(1-\delta_{m,n-1}) q_n + t^{m-1}(1-\delta_{m,1}) q_n.
\end{split}
\end{equation*}
Hence
\begin{equation}\label{eq:action}
T(t^m q_n) = \sum_{i=1}^{n-1} \alpha_{m,i} t^iq_n.
\end{equation}

For $1 \le m \le n-1$, define
$r_{m,n}:= \sum_{i=1}^{n-1} \beta_{m,i}^{(n)} t^i q_n$ and
$p_{m,n}:= r_{m,n} r_{m,n}^*$.
Observe that, if we identify the $i$th standard basis vector with
$t^iq_n$, $1 \le i \le n-1$, then $r_{m,n}$ is an eigenvector of $A_n$ with
eigenvalue $\lambda_{m,n}$.  Moreover, we have just
checked that $T$ acts like $A_n$ on the span of the $t^m q_n$. This partially
explains why the $r_{m,n}$ give rise to pairwise orthogonal projections
with image contained in the eigenspace of $T$ for the eigenvalue
$\lambda_{m,n}$, which is essentially the statement of the following
lemma.

\begin{lemma}\label{lem:key}
$(p_{m,n} \mid n \ge 2, 1 \le m \le n -1)$ is a family of pairwise
orthogonal projections in $\complexs[U\wr \integers]$ which is complete,
that is,
$\sum_{n \ge 2} \sum_{m=1}^{n-1} \tr_{U \wr \integers} (p_{m,n}) = 1.$
Moreover, if $1 \le m\le n-1$, then $T(p_{m,n}) = \lambda_{m,n}p_{m,n}$.
\end{lemma}

\begin{proof} Let $1 \le m \le n-1$ and $1 \le m' \le n'-1$.

Here $$r_{m,n}^* =  q_n^* \sum_{i=1}^{n-1} (t^i)^*\beta_{m,i}^{(n) *}
= q_n \sum_{i=1}^{n-1}  t^{-i} \beta_{m,i}^{(n)}.$$
Thus
\begin{equation*}
\begin{split}
r_{m',n'}^* r_{m,n} &=  q_{n'} \sum_{j=1}^{n'-1} t^{-j} \beta_{m',j}^{(n')}
\sum_{i=1}^{n-1} \beta_{m,i}^{(n)} t^i q_n \\
& =  \delta_{n,n'} q_n \sum_{i=1}^{n-1} \beta_{m',i}^{(n)} \beta_{m,i}^{(n)}
\text { by Lemma~\ref{lem:orthogonality}}\\ & = \delta_{n,n'} q_{n}
\delta_{m,m'}\text { by~(\ref{eq:orthogonal}).}
\end{split}
\end{equation*}

It follows that the $p_{m,n}$ are pairwise orthogonal.

Moreover,
\begin{equation*}
\begin{split}
\tr_{U \wr \integers} (p_{m,n})  &=
\tr_{U \wr \integers} (r_{m,n} r_{m,n}^*) =
\tr_{U \wr \integers} (r_{m,n}^* r_{m,n}) \\&=
\tr_{U \wr \integers} (q_n)
= \frac{(W-1)^2}{W^n} \text{ by~(\ref{eq:trace})}.
\end{split}
\end{equation*}

Now,
\begin{equation*}
\begin{split}
\sum_{n \ge 2} \sum_{m=1}^{n-1}  \tr_{U \wr \integers} (p_{m,n}) &=
\sum_{n\ge 2} \sum_{m=1}^{n-1} \frac{(W-1)^2}{W^n}
=\sum_{n\ge 2} (n-1)\frac{(W-1)^2}{W^n} \\
&=\sum_{n\ge 1} n \frac{(W-1)^2}{W^{n+1}}
= (1 - \frac 1 W)^2 \sum_{n\ge 1} n (\frac{1}{W})^{n-1} = 1,
\end{split}
\end{equation*}
since, for $\abs{x}<1$, $\sum_{n\ge 1} n x^{n-1}= (\sum_{n\ge 0}
x^n)' = (\frac{1}{1-x})' = \frac{1}{(1-x)^2}$.

Also,
\begin{equation*}
\begin{split}
T(r_{m,n}) &=  T(\sum_{j=1}^{n-1} \beta_{m,j}^{(n)}t^jq_n)
= \sum_{j=1}^{n-1} \beta_{m,j}^{(n)} T(t^jq_n) \\
&= \sum_{j=1}^{n-1} \beta_{m,j}^{(n)} \sum_{k=1}^{n-1} \alpha_{j,k} t^kq_n
\text{
by~(\ref{eq:action}) } \\ & = \sum_{k=1}^{n-1}
(\sum_{j=1}^{n-1}\beta_{m,j}^{(n)}
\alpha_{j,k}) t^kq_n
\\ & =  \sum_{k=1}^{n-1} \lambda_{m,n} \beta_{m,k}^{(n)}  t^kq_n \text{
by~(\ref{eq:diagonal})}\\ & =  \lambda_{m,n} r_{m,n}.
\end{split}
\end{equation*}
Thus $T(r_{m,n}) = \lambda_{m,n} r_{m,n}$, and, on right multiplying by
$r_{m,n}^*$, we see $T(p_{m,n}) = \lambda_{m,n} p_{m,n}$.
\end{proof}

We have now `diagonalized' $T$ in the sense that we have decomposed
$l^2(U \wr \integers)$ into the Hilbert sum of subspaces of the form
$p_{m,n}(l^2(U\wr\integers))$ on which $T$ acts as multiplication by the scalar
$\lambda_{m,n}$.

Hence, for each $\mu \in \reals$, $\ker(T - \mu)$ is the Hilbert sum of those
$p_{m,n}(l^2(U \wr \integers))$ such that $\lambda_{m,n} = \mu$.  Thus
either $\ker(T - \mu) = 0$  or $\mu = \lambda_{m_0,n_0}$ for some  $m_0$,
$n_0$ with
$1 \le m_0 \le n_0-1$.

We now consider the latter case.  Here, for all $(m,n)$, $\lambda_{m,n} =
\mu$ if and
only if $\frac m n = \frac {m_0} {n_0}$.  We may assume that $m_0$ and $n_0$ are
coprime, so $\mu \in M_{n_0}$.  Also,  $\lambda_{m,n} = \mu$ if and only if
$(m,n) = (im_0,in_0)$ for some $i \ge 1$.  Thus $\ker(T - \mu)$ is the
Hilbert sum of
the $p_{im_0,in_0}(l^2(U \wr\integers))$ with $i \ge 1$; hence
\begin{equation*}
\begin{split}
\dim_{U\wr \integers}(\ker(T- \lambda_{m_0,n_0})) &=
\sum_{i \ge 1} \dim_{U\wr \integers}(p_{im_0,in_0}(l^2(U \wr \integers)))
\\=\sum_{i \ge 1} \tr_{U\wr \integers}(p_{im_0,in_0})& =
\sum_{i \ge 1}  \frac{(W-1)^2}{W^{in_0}}=
\frac{(W-1)^2}{W^{n_0} -1}.
\end{split}
\end{equation*}

Theorem~\ref{bigmain} now follows.

\begin{remarks}
   The hypothesis in Theorem~\ref{bigmain} that $U$ has torsion could be
   weakened to the assumption that $\complexs[U]$ has a nontrivial projection;
   however, if $U$ is torsion-free, it is conjectured, and known in many cases,
   that $\complexs[U]$ does not contain any nontrivial projections.

   It easy to show that the hypothesis in Theorem~\ref{bigmain} that $e$ is a
   nontrivial projection in $\complexs[U]$ can be weakened to the assumption that $e$
   is a nontrivial projection in $\NeumannN(U)$; here, the hypothesis that
$U$ has
   torsion should be weakened to the assumption that $U$ is nontrivial. \qed
\end{remarks}

\section{Direct products of wreath products}

We now produce even more unusual examples by taking direct products of the
groups
studied so far.

\begin{theorem}\label{smallmain}
  Let $U$ and $V$ be groups with torsion, and
$G = (U \wr \integers) \times (V \wr \integers)$.  Let $e$ be a nontrivial
projection
in $\complexs [U]$ and $f$ a nontrivial projection in $\complexs [V]$.   Let
$X = (\tr_U(e))^{-1}$ and  $Y = (\tr_V(f))^{-1}$, so $X > 1$, $Y >1$.  Let
$T = T(U,e) \in \complexs[U \wr \integers] \subset \complexs[G]$, and
$S = T(V,f) \in \complexs[V\wr\integers]\subset \complexs[G]$.
Then
\begin{equation}\label{double_sum}
\begin{split}
&\dim_G(\ker(T-S)) \\&=
(X-1)^2(Y-1)^2(\sum_{m\ge 1}\sum_{n \ge 1}\frac{\gcd(m,n)}{X^mY^n}) -
(X-1)(Y-1).
\end{split}
\end{equation}
\end{theorem}

\begin{proof} By Lemma~\ref{lem:key}, there is a complete family
$(p_{m,n} \mid n \ge 2, 1 \le m < n)$ of pairwise
orthogonal projections in $\complexs[U\wr \integers]$, such that,
if $1 \le m < n$, then $T(p_{m,n}) = \lambda_{m,n}p_{m,n}$, and,
by~(\ref{eq:trace}),
$\tr_{U \wr \integers}(p_{m,n}) = \frac{(X-1)^2}{X^n}$.

Similarly, there is a complete family
$(q_{m,n} \mid n \ge 2, 1 \le m < n)$ of pairwise
orthogonal projections in $\complexs[V\wr \integers]$ such that, if $1 \le
m < n$,
then $S(q_{m,n}) = \lambda_{m,n}q_{m,n}$, and
$\tr_{V \wr \integers}(q_{m,n}) = \frac{(Y-1)^2}{Y^{n}}$.

By Lemma~\ref{lem:tensor_product_trace}, there is a complete family
$$(p_{m,n}q_{m',n'} \mid n,n' \ge 2, 1 \le m < n, 1 \le m' < n')$$ of pairwise
orthogonal projections in $\complexs[G]$, such that,
if $1 \le m < n$ and $1 \le m' < n'$ then
$$ T(p_{m,n}q_{m',n'}) = \lambda_{m,n}p_{m,n}q_{m',n'}\text{ and }
S(p_{m,n}q_{m',n'}) = \lambda_{m',n'}p_{m,n}q_{m',n'},$$ and
$$\tr_G(p_{m,n}q_{m',n'})
=  \frac{(X-1)^2}{X^n}\frac{(Y-1)^2}{Y^{n'}}.$$

Thus $l^2(G)$ is the Hilbert sum of the subspaces of the form
$p_{m,n}q_{m',n'}(l^2(G))$ where $T-S$ acts as
multiplication by the scalar $\lambda_{m,n} - \lambda_{m',n'}$.

Hence $\ker(T-S)$ is the Hilbert sum of the
$p_{m,n}q_{m',n'}(l^2(G))$ such that $\lambda_{m,n} = \lambda_{m',n'}$.

Therefore, $$\dim_G(\ker(T-S))
= \sum_{n\ge 1} \sum_{n'\ge 1} b(n,n')
\frac{(X-1)^2}{X^n}\frac{(Y-1)^2}{Y^{n'}}$$
where $b(n,n')$ is the number of pairs $(m,m')$ such that
$1 \le m < n$, $1 \le m' <n'$, and $\frac{m}{n} = \frac{m'}{n'}$.  But such
pairs
correspond bijectively to the fractions of the form
$\frac{m_0}{\gcd(n,n')}$, $ 1 \le m_0 < \gcd(n,n')$.  Thus
$b(n,n') = \gcd(n,n') - 1$.  Hence
\begin{equation*}
\begin{split}
&\dim_G(\ker(T-S))
= \sum_{n\ge 1} \sum_{n'\ge 1}
\frac{(\gcd(n,n')-1)(X-1)^2(Y-1)^2}{X^{n}Y^{n'}}\\
&=\sum_{n\ge 1} \sum_{n'\ge 1}
\frac{\gcd(n,n')(X-1)^2(Y-1)^2}{X^{n}Y^{n'}} -
\sum_{n\ge 1} \sum_{n'\ge 1}
\frac{(X-1)^2(Y-1)^2}{X^{n}Y^{n'}}.
\end{split}
\end{equation*}

Since $\sum_{n\ge 1} \frac{1}{X^n} = X^{-1}\frac{1}{1-X^{-1}} = \frac{1}{X
- 1}$,
the result follows.
\end{proof}

\begin{remarks}\label{gcd_series}  Recall that, for any positive integer $n$,
$\phi(n)$ denotes the number of primitive $n$th roots of unity, so
$\abs{M_n} =\phi(n)$.

For $X > 1$, $Y >1$, the double infinite sum
occurring in (\ref{double_sum}) has an expession as a single infinite sum,
$$\sum_{m\ge 1}\sum_{n \ge 1}\frac{\gcd(m,n)}{X^mY^n}
=\sum_{k\ge 1} \frac{\phi(k)}{(X^k-1)(Y^k-1)},$$
since
\begin{equation*}
\sum_{k\ge 1}\frac{\phi(k)}{(X^k-1)(Y^k-1)} =
\sum_{k\ge 1}\phi(k)\sum_{i\ge 1} X^{-ik} \sum_{j\ge 1} Y^{-jk}
= \sum_{m\ge 1} \sum_{n\ge 1}  \frac{a(m,n)}{X^mY^n}
\end{equation*}
where
\begin{equation*}
a(m,n) = \sum_{\{k \ge 1 :  k \vert m, k \vert n\}} \phi(k)
= \sum_{k\vert\gcd(m,n)} \phi(k) = \gcd(m,n).
\end{equation*}
It follows that $$\dim_G(\ker(T-S)) =
(X-1)^2(Y-1)^2 \sum_{k\ge 2} \frac{\phi(k)}{(X^k-1)(Y^k-1)}. \qed $$

\end{remarks}

\section{$L^2$-Betti numbers}
\label{sec:l2-betti-numbers}

We previously observed that, by results of Kaplansky and Zaleskii, the traces of
projections in complex, or rational, group algebras are rational numbers in the
interval $[0,1]$. In order to maximize the scope of Theorem~\ref{smallmain} for
producing examples of $L^2$-Betti numbers, we need the following result which
shows that the traces of projections in rational group algebras are
{\it precisely} the rational numbers in the interval $[0,1]$.   We write
$C_n$ for a cyclic group of order $n$, written multiplicatively, with generator
$t = t_n$.

\begin{lemma}\label{existence_of_rational_projections}
  Let $q$ be a rational number in the interval $[0,1]$.  Then there is an
expression
  $q = \frac{m}{n}$ where the denominator has the form
  $n = 2^rs$ with $s$ odd and $2^r \ge s-1$, and, for any such expression,
  $\rationals[C_n]$ contains some projection $e$ with trace $q$, and $ne \in
  \integers[C_n]$.
\end{lemma}

\begin{proof}
  By multiplying the numerator and denominator of $q$ by a sufficiently
high power of
  $2$, we see that $q$ has an expression of the desired type.  Now consider any
  expression $q = \frac{m}{n}$ where $n = 2^rs$ with $s$ odd and $2^r \ge s-1$.

  We first show, by induction on $r$, that, if $0 \le c \le 2^{r}$, then
  $\rationals[C_{2^r}] = \rationals[t \mid t^{2^{r}}=1]$
  has an ideal whose dimension over $\rationals$ is $c$.  Since the
  orthogonal complement is then an ideal of dimension $2^r - c$ over the
  rationals, it amounts to the same if we consider only $c \le 2^{r-1}$.
For $r = 0$,
  we can take the zero ideal;  thus,  we may assume that $r\ge 1$ and the
result holds
  for smaller $r$.  Now $\rationals[C_{2^r}]$ has a projection
  $e =\frac{1+t^{2^{r-1}}}{2}$;  this is $\avg(U)$ for the subgroup $U$ of
order $2$
  in $C_{2^r}$.  As rings
  $$e\rationals[C_{2^r}] \simeq \rationals[C_{2^r}]/(1-e)
  \simeq \rationals[C_{2^{r-1}}].$$  By the induction hypothesis, the
latter has an
  ideal of dimension $c$ over $\rationals$, and viewed in $e\rationals[C_{2^r}]$
  this is an ideal of $\rationals[C_{2^r}]$.  This completes the proof by
  induction.  Hence, if $0 \le c \le 2^{r}$, then $\rationals[C_{2^r}]$ has a
  projection  $e(c)$ with $\tr_{C_{2^r}}(e(c)) = \frac{c}{2^r}$.

  Let $f = \avg(C_s) \in \rationals[C_s]$, so $\tr_{C_s}(f) = \frac{1}{s}$, and
  $\tr_{C_s}(1-f) = \frac{s-1}{s}$.

  By identifying
  \begin{equation*}
\rationals[C_n] = \rationals[C_n^s \times C_n^{2^r}]
  = \rationals[C_{2^r} \times C_s],
\end{equation*}
we see that, for $0 \le c \le 2^r$, we have
  projections $e(c)f$ and $e(c)(1-f)$ in $\rationals[C_n]$, with traces
  $\frac{c}{2^r}\frac{1}{s} = \frac{c}{n}$ and
  $\frac{c}{2^r}\frac{s-1}{s} = \frac{c(s-1)}{n}$, respectively, by
  Lemma~\ref{lem:tensor_product_trace}.

  We claim there exist integers $a$, $b$ with $0 \le a,b \le 2^r$ such that
  $a + (s-1)b = m$.  We know that $0 \le m \le n = 2^rs$.  If $m \ge 2^r(s-1)$,
  then $m \in [2^r(s-1), 2^rs]$, and we can take $b = 2^r$ and
  $a = m - (s-1)b = m - 2^r(s - 1) \in [0,2^r]$. If $m < 2^r(s-1)$, then, by
  the division algorithm,  $m = (s-1)b + a$ with $0  \le b < 2^r$, and
  $0 \le a \le s-2 < 2^r$.  This proves the claim.

  Now let $e = e(a)f + e(b)(1-f)$, a sum of orthogonal projections.  Thus,
$e$ is
  a projection and
  $$\tr_{C_n}(e) =  \tr_{C_n}(e(a)f) + \tr_{C_n}(e(b)(1-f))
  = \frac{a}{n} + \frac{b(s-1)}{n} = \frac{a + b(s-1)}{n} = \frac{m}{n},$$
  as desired.

  It remains to show that $e$ lies in  $\frac{1}{n}\integers[C_n]$, but it
is well
  known that this holds for all the idempotents of $\rationals[C_n]$.
  Alternatively, it is straightforward to check that all the projections
involved in
  the foregoing proof have the right denominators.
\end{proof}

We now obtain the following special case of Theorem~\ref{smallmain}.

\begin{corollary}\label{cor:products}
Let $p$ and $q$ be rational numbers with $0 < p, q < 1$.  There exist positive
integers $m$ and $n$, and projections
\begin{equation*}
e=e^*=e^2 \in\rationals[C_m], \qquad f
=f^*=f^2\in\rationals[C_n]
\end{equation*}
with $\tr_U(e)=p$, $\tr_V(f)=q$.  Let
 \begin{equation*}
G(p,q) := (C_m \wr \integers) \times (C_n \wr \integers),
\end{equation*}
\begin{equation*}
T: = T(U,e) \in \complexs[U \wr \integers] \subset \complexs[G],\text{
  and }
S := T(V,f) \in \complexs[V\wr\integers]\subset \complexs[G].
\end{equation*}
Let $Z = Z(p,q):= mn(T-S)$, and let
\begin{equation*}
\begin{split}
\kappa & = \kappa(p,q):=(p^{-1}-1)^2(q^{-1}-1)^2 \sum_{k\ge
2}\frac{\phi(k)}{(p^{-k}-1)(q^{-k}-1)} \\&=
(p^{-1}-1)^2(q^{-1}-1)^2 \left(\sum_{i\ge 1}\sum_{j\ge 1}
\gcd(i,j) p^{i} q^{j}\right) - (p^{-1}-1)(q^{-1}-1).
\end{split}
\end{equation*}
Then $Z \in \integers[G]$ and $\dim_G(\ker Z) = \kappa$. \qed
\end{corollary}

\begin{remarks}\label{Betti_numbers}
Let $0<p,q<1$ be rational numbers.  Let $G = G(p,q)$, $Z = Z(p,q)$ and
$\kappa = \kappa(p,q)$ as in Corollary~\ref{cor:products}.

By the Higman Embedding Theorem, any recursively presented group can be
embedded in a
finitely presented group, so $G$ can be embedded in a finitely presented group
$H$.  (Here it is easy to find an explicit suitable finitely presented
group; see, for
example,~\cite{Baumslag(1972)}
or~\cite[Lemma 3]{Grigorchuk-Linnell-Schick-Zuk(2000)}. This explicit
supergroup has the additional nice property of being metabelian, that is,
$2$-step
solvable.  Moreover, one can precisely describe its finite subgroups.)

By Corollary~\ref{cor:products}, $Z \in\integers[G]\subseteq \integers[H]$ and
$\dim_H(\ker Z) =  \dim_G(\ker Z) = \kappa$.

It is then well known how to construct a finite CW-complex or a closed
manifold $M$
with $\pi_1(M)\simeq H$ and with third $L^2$-Betti number $\kappa$;  see, for
example,~\cite{Grigorchuk-Linnell-Schick-Zuk(2000)}.

Thus $\kappa(p,q)$ is an $L^2$-Betti number of a closed manifold.  It is
conceivable that
this is a counterexample to Atiyah's conjecture~\cite{Atiyah(1976)} that
$L^2$-Betti
numbers of closed manifolds are  rational, but we have not been able to
decide whether
$\kappa(p,q)$ is rational or not.\qed
\end{remarks}

\begin{example}
Consider $\kappa(\frac 12, \frac 12) = \sum_{k\ge2}\frac{\phi(k)}{(2^k-1)^2} =
0.1659457149\ldots\, .$ If we sum the first 400 terms,  then elementary
methods show
that the remaining tail is less than $10^{-201}$.  This allows us to
calculate the
first 199 terms of the continued fraction expansion of $\kappa(\frac 12,
\frac 12)$.
One consequence we find is that if $\kappa(\frac 12, \frac 12)$ is rational
then both
the numerator and
the denominator exceed $10^{100}$.  It seems reasonable to assert that
$\kappa(\frac 12, \frac 12)$
is not obviously rational.\qed
\end{example}

\section{Power series}

Throughout this section, let $\complexs((x,y))$ denote the field of
(formal) Laurent series in two variables (with complex coefficients).

The expression
\begin{equation*}
\Phi(x,y):=  \sum\limits_{m\ge 1}\sum\limits_{n\ge 1}\gcd(m,n)x^my^n
\end{equation*}
arising from (\ref{double_sum}) can be viewed as an
element of $\complexs((x,y))$.  By Remarks~\ref{Betti_numbers}, if there exist
rational numbers $p$, $q$ in the interval $(0,1)$ such that (the limit of)
$\Phi(p,q)$  is irrational, then there exists a counterexample to the Atiyah
conjecture; so it is of interest to know whether $\Phi(p, q)$ is always
rational for
such rational numbers $p,q$.  One
(traditionally successful) way to show that such an expression is rational
would be
to show that $\Phi(x, y)$ itself is rational, that is, lies in the
subfield $\rationals(x,y)$ of
rational Laurent series over the rationals. In this section, we will
eliminate this
possibility by showing that  $\Phi(x, y)$ is transcendental over
$\complexs(x,y)$.  In fact, we will show the stronger result that the
specialization
$\Phi(x, x)$ is transcendental over $\complexs(x)$.

The following result is well known, but we have not found a reference.  The
proof is left to the reader.

\begin{lemma}\label{diffeq}
Suppose that $f \in \complexs((x))$ is algebraic over $\complexs(x)$ of
degree $d$.
Then the subfield $\complexs(x,f)$ is closed under the usual derivation
operation,
$F\mapsto F' = \frac{dF}{dx}$, on $\complexs((x))$.  Moreover,
$\complexs(x,f)$ is a
$d$-dimensional vector space over $\complexs(x)$, so the $d+1$ higher-order
derivatives $f^{(i)}:= (\frac{d}{dx})^i(f)$,  $0 \le i \le d$, are
$\complexs(x)$-linearly dependent.  Hence
$f$ satisfies some non-trivial order $d$ differential equation over
$\complexs(x)$. \qed
\end{lemma}

We can now apply this lemma to get a transcendentality criterion.

\begin{proposition}\label{gaps}
Suppose that $a\colon \naturals \to \complexs$, $n \mapsto a(n)$, has the
property
that, for each $N \in \naturals$, there exist infinitely many $m \in
\naturals$ such
that, whenever $j \in \integers$ satisfies $1 \le \abs{j} \le N$,
$$\abs{a(m)} > N
\abs{a(m+j)}.$$  Then the power series
$\sum\limits_{n \ge 0} a(n) x^n \in \complexs((x))$
does not satisfy any non-trivial differential equation over
$\complexs(x)$, so is
transcendental over $\complexs(x)$.
\end{proposition}

\begin{proof} Let $f:= \sum_{n \ge 0} a(n) x^n \in \complexs((x))$, and
suppose that
$f$ satisfies a non-trivial differential equation over $\complexs(x)$,
\begin{equation}\label{eq:diff}
\sum_{i=0}^d q_i f^{(i)} = 0
\end{equation}
where $q_i \in \complexs(x)$, not all zero.  By multiplying through by a common
denominator, we may assume that all the $q_i$ lie in
$\complexs[x]$.  (Notice it is natural not to have a ``constant term"
on the right-hand side of (\ref{eq:diff}) since it could
be eliminated by iterated derivation of the equation.)

Viewing~(\ref{eq:diff}) as a collection of equations, one for each power
$x^n$, we see
that there exists some $N \in \naturals$, and polynomials $p_k(t) \in
\complexs[t]$
such that
\begin{equation}\label{relations}
\sum_{k=0}^N p_k(n)a(n+k) = 0 \text{ for all } n\in \naturals.
\end{equation}

Choose $k_0$, with $0 \le k_0 \le N$, and $n_0 \in \naturals$ such that
$\abs{p_{k_0}(n)} \ge \abs{p_k(n)}$
for all
$n \ge n_0$, and all $k$ with $0 \le k \le N$. In other words, $p_{k_0}$
eventually
dominates all
the $p_k$, $0 \le k \le N$.

It follows from the hypothesis on the $a(n)$ that there exists
 $m \in \naturals$ such that $m \ge n_0 + k_0$, and $\abs{a(m)} > N
\abs{a(m+j)}$ for
all  $j \in \integers$ with  $1 \le \abs{j} \le N$.  Now take $n = m - k_0$.
Then $n \ge n_0$, and
$$\abs{a(n+k_0)} > \sum_{k=0}^{k_0-1} \abs{a(n+k)} +
\sum_{k=k_0+1}^{N} \abs{a(n+k)}.$$  Thus
\begin{equation*}
\begin{split}
\abs{p_{k_0}(n)a(n+k_0)} &> \sum_{k=0}^{k_0-1} \abs{p_{k}(n)a(n+k)} +
\sum_{k=k_0+1}^{N} \abs{p_{k}(n)a(n+k)}\\
&\ge \abs{(\sum_{k=0}^{N} p_{k}(n)a(n+k)) - p_{k_0}(n)a(n+k_0)} \\&=
\abs{0-p_{k_0}(n)a(n+k_0)} \text{ by~(\ref{relations}).}
\end{split}
\end{equation*}

This contradiction shows that $f$ does not satisfy any non-trivial differential
equation over
$\complexs(x)$, so, by Lemma~\ref{diffeq}, $f$ is not algebraic over
$\complexs(x)$.
\end{proof}

We now record some important results from number theory that we shall require.

\begin{lemma}\label{largeQ}
For each positive integer $i$, let $p_i$ denote the $i$th prime
number.  There
exists an integer $Q_0$ such that, for all $Q \ge Q_0$, the following hold.
\begin{enumerate}
\item\label{factorial}  $Q! \le (\frac Q2)^Q.$
\item\label{prime_number_theorem} $\frac 3 4 \le \frac{p_Q}{Q \log Q} \le
\frac 54$.
 \item\label{Mertens} $\prod\limits_{i=1}^Q ( 1 - \frac{1}{p_i}) \ge \frac 1 Q$.
\end{enumerate}
\end{lemma}

\begin{proof}
In the following, $f(Q)=o(g(Q))$ means $\lim_{Q\to\infty} f(Q)/g(Q)
=0$, and $f(Q)\sim g(Q)$ means $\lim_{Q\to\infty} f(Q)/g(Q) =1$.

\ref{factorial}  By Stirling's formula,
 $Q! \sim \sqrt{2\pi}Q^{Q+\frac12}e^{-Q}$, and the latter is $o((\frac
Q2)^Q)$, since
$e > 2$.  One can argue directly that $\sum_{i=1}^Q \log i \le
\int_1^{Q+1} \log
x~dx$, so  $$\log Q! \le (Q+1)\log (Q+1) - Q,$$ so
$$Q! \le (Q+1)^{Q+1}e^{-Q}= Q^Q(1+\frac{1}{Q})^Q(Q+1)e^{-Q} = o((\frac
Q2)^Q),$$
since $e > 2$.

\ref{prime_number_theorem}  By the Prime Number Theorem, $p_Q \sim Q\log Q$;
see~\cite[Theorem~8, pages~10,~367]{Hardy-Wright(1979)}.

\ref{Mertens} By Mertens' Theorem,
$\prod\limits_{i=1}^Q ( 1 - \frac{1}{p_i}) \sim \frac {e^{-\gamma}}{\log
p_Q}$, where
$\gamma$ is Euler's constant; see~\cite[Theorem~429,
page~351]{Hardy-Wright(1979)}.
By the Prime Number Theorem, $\log p_Q \sim \log Q$, so
$\prod\limits_{i=1}^Q ( 1 - \frac{1}{p_i}) \sim \frac {e^{-\gamma}}{\log
Q}$.  Since
$\frac 1Q = o(\frac{1}{\log Q})$,  we see that
$\frac 1Q = o(\prod\limits_{i=1}^Q ( 1 - \frac{1}{p_i}))$.

The result now follows.\end{proof}

\begin{theorem}
$\Phi(x,x) = \sum\limits_{m\ge 1} \sum\limits_{n\ge 1} \gcd(m,n) x^{m+n}$ and
$\sum\limits_{n \ge 1} \sum\limits_{d\vert n}\frac{\phi(d)}{d} nx^n$ are
\linebreak transcendental over
$\complexs(x)$.
\end{theorem}

\begin{proof}  For each positive integer $n$, let
$a(n):= n \sum_{d\vert n} \frac{\phi(d)}{d}$.
Thus
\begin{equation*}
\begin{split}
a(n)&= n \sum_{d\vert n} \frac{\phi(d)}{d} = \sum_{d\vert n} \frac{n}{d}
\phi(d)
= \sum_{d\vert n}  \sum_{\{i : 1 \le i \le n, d \vert i\}} \phi(d)
= \sum_{i=1}^n \sum_{\{d : d\vert i, d\vert n\}} \phi(d) \\
& = \sum_{i=1}^n \gcd(i,n) =  \sum_{i=1}^n \gcd(i,n-i)
= \sum_{i=1}^{n-1} \gcd(i,n-i) + n.
\end{split}
\end{equation*}
Now
\begin{equation*}
\Phi(x,x) = \sum_{i \ge 1}\sum_{j\ge 1} \gcd(i,j) x^{i+j}
= \sum_{n \ge 1}\sum_{i=1}^{n-1} \gcd(i,n-i) x^{n},
\end{equation*}
so
\begin{equation*}
(\sum_{n \ge 1} a(n) x^{n}) - \Phi(x,x) =
\sum_{n \ge 1} n x^{n} =  x(\sum_{n \ge 0} x^{n})' =  \frac{x}{(1-x)^2}.
\end{equation*}
Thus $\sum_{n \ge 1} a(n) x^{n}$ and  $\Phi(x,x)$ differ by an element of
$\rationals(x)$, so it suffices to show that
$\sum_{n\ge 1} a(n) x^n$ is transcendental over $\complexs(x)$.

By Proposition~\ref{gaps}, it suffices to show that, for each
$N \in \naturals$, there exist infinitely many $m \in \naturals$ such
that, whenever
$j \in \integers$ satisfies $1 \le \abs{j} \le N$,  $$\abs{a(m)} > N
\abs{a(m+j)}.$$

We may suppose that $N$ is fixed.

Remember the $p_i$ is the $i$th prime number. For each $Q\in\naturals$,
let
\begin{equation*}
m_Q:= \prod_{i=1}^Q p_i \prod_{i=1}^N p_i^N.
\end{equation*}

We may now suppose that $j$ is fixed with $1 \le \abs{j} \le N$, and it
suffices to
show that
\begin{equation*}
\lim\limits_{Q\to\infty} \frac{a(m_Q+j)}{a(m_Q)} = 0.
\end{equation*}

We use the notation of Lemma~\ref{largeQ}, concerning $Q_0$.
Let
\begin{equation*}
 C_1 =  \prod_{i=1}^{Q_0} p_i \prod_{i=1}^N p_i^N.
\end{equation*}

Now suppose that $Q$ is an integer with  $Q \ge \max\{Q_0,N\}$, let $m =
m_Q$ and let
$m' = \prod_{i=1}^Q p_i.$

We wish to bound $a(m) = m \sum_{d \vert m} \frac{\phi(d)}{d}$ from below.
Recall
that, for  any positive integer $n$,  $\frac{\phi(n)}{n} = \prod (1 - \frac
1p)$,
where the product is over the distinct prime divisors $p$ of
$n$.  Thus $a(m) \ge m \sum_{d \vert m} \frac{\phi(m)}{m} = m
\de(m)\frac{\phi(m)}{m}$,
where $\de(m)$ denotes the number of divisors $d$ of $m$.  Also,
$\frac{\phi(m)}{m} = \prod\limits_{i=1}^Q ( 1 - \frac{1}{p_i})$, which, by
Lemma~\ref{largeQ}\ref{Mertens}, is at least~$\frac 1 Q$.   Thus
$a(m) \ge m \de(m) \frac 1Q$.  Notice that $\de(m) \ge \de(m')$, since $m'$
divides $m$.
From the definition of $m'$, we see that $\de(m') = 2^Q$.  Thus
$$a(m) \ge m 2^Q \frac 1Q.$$

We next wish to bound $a(m+j)$ from above.  Let $\Omega(m+j)$ be the
number, counting multiplicity, of prime factors of $m$, and
let
$$m+j = p_{i_1}p_{i_2}\cdots p_{i_{\Omega(m+j)}}$$
be the factorization of $m+j$ into prime factors. Then
$\de(m+j) \le 2^{\Omega(m+j)}$, and
\begin{equation*}
\begin{split}
a(m+j) &= (m+j)\sum_{d \vert (m+j)} \frac{\phi(d)}{d}
\le  (m+j)\sum_{d\vert (m+j)} 1 = (m+j) \de(m+j)
\\&\le (m+j) 2^{\Omega(m+j)} \le (m+N)2^{\Omega(m+j)} \le 2m 2^{\Omega(m+j)}.
\end{split}
\end{equation*}
Consider $ 1 \le l \le \Omega(m+j)$.  If $i_l \le Q$, then $p_{i_l}$
divides $m$ so
$p_{i_l}$ divides $j$.  But $1\le \abs{j} \le N$, so $p_{i_l} \le N$, so
$i_l \le N$.  Hence $p_{i_l}^N$ divides $m$, but $p_{i_l}^N \ge 2^N > N \ge
\abs{j}$,
so $p_{i_l}^N$ cannot divide $j$, so cannot divide $m+j$.  Thus, the number of
$i_l$ which are less than $Q$ is at most $N^N$.  Let $z = z(Q,j)$ denote
the number of
$l$ such that $i_l \ge Q$, so $\Omega(m+j) \le z+ N^N,$  and
$$a(m+j) \le 2m 2^{\Omega(m+j)} \le 2m 2^{z+N^N}.$$

Thus
$$\frac{a(m+j)}{a(m)} \le \frac{2m 2^{z+N^N}}{m2^Q\frac{1}{Q}} = Q2^{z-Q}
2^{N^N+1}.$$
Hence it remains to show that $\lim\limits_{Q\to \infty} Q2^{z-Q}=0$, or
equivalently, $$\lim\limits_{Q\to \infty} Q - z - \log_2Q = \infty.$$

Since $z$ is the number, counting multiplicity, of prime factors $p_{i_l}$
of $m+j$
with $p_{i_l}\ge p_Q$,
$$p_Q^z \le m + j \le m + N \le 2m.$$ We can write
\begin{equation*}
\begin{split}
m &= \prod_{i=1}^Q p_i \prod_{i=1}^N p_i^N
\le \prod_{i=1}^{Q_0} p_i \prod_{i=2}^Q (\frac {5}{4} ~ i \log i )
\prod_{i=1}^N
p_i^N = C_1 \prod_{i=2}^Q (\frac {5}{4} ~ i \log i ) \\
&\le C_1 (\frac 54)^Q Q!(\log Q)^Q
\le C_1 (\frac 54)^Q (\frac Q2)^Q (\log Q)^Q,
\end{split}
\end{equation*}
by Lemma~\ref{largeQ}\ref{factorial}.
Thus $$(\frac{3}{4}~Q \log Q)^z \le p_Q^z \le 2m
\le 2C_1 (\frac 54)^Q (\frac Q2)^Q (\log Q)^Q.$$
Hence $$(\frac 34 Q \log Q)^{z-Q}
\le 2C_1 (\frac 43)^Q (\frac 12)^Q (\frac 54)^Q = 2C_1 (\frac 56)^Q,$$
so $(z-Q)(\log \frac 34 + \log Q + \log \log Q) \le \log 2C_1 - Q \log(
\frac 65)$,
and  $$-(Q-z) \le
\frac{\log 2C_1 - Q \log( \frac 65)}{\log \frac 34 + \log Q + \log \log Q}
 \sim -\log( \frac 65) \frac{Q}{\log Q}.$$
It follows that
$$\lim\limits_{Q\to \infty} Q - z - \log_2Q \ge
\lim\limits_{Q\to \infty}\log( \frac 65) \frac{Q}{\log Q} -  \log_2Q  =
\infty,$$
as desired.
\end{proof}

\bibliographystyle{plain}
\bibliography{DicksSchick}

\end{document}